\documentclass{amsart}
\usepackage{amssymb,amsmath, amsthm,latexsym}
\usepackage{graphics}
\usepackage{amscd}
\usepackage{graphics}
\usepackage{here}
\newcommand{\cal}[1]{\mathcal{#1}}
\theoremstyle{plain}

\newtheorem{lemma}{Lemma}[section]

\newtheorem{proposition}[lemma]{Proposition}

\parskip=\bigskipamount

\let\egthree=\phi
\let\phi=\varphi
\let\varphi=\egthree

\begin{document}
\title
{Closed Teichm\"uller geodesics in the thin part of 
moduli space} 
\author{Ursula Hamenst\"adt}
\thanks
{
Research
partially supported by Sonderforschungsbereich 611}
\date{August 8, 2009}

\begin{abstract}
Let $S$ be an oriented surface of 
genus $g\geq 0$ with $m\geq 0$ punctures
where $3g-3+m\geq 2$.
We show that if 
$3g-3+m\geq 4$
then for
every compact subset $K$ of the moduli
space ${\rm Mod}(S)$ for $S$ there is a closed
Teichm\"uller geodesic in ${\rm Mod}(S)$
which does not intersect $K$.
\end{abstract}

\maketitle

\section{Introduction}

Let $S$ be a \emph{nonexceptional} 
oriented surface of finite type, i.e. 
$S$ is a
closed surface of genus $g\geq 0$ from which $m\geq 0$
points, so-called \emph{punctures},
have been deleted, and $3g-3+m\geq 2$.
The \emph{Teichm\"uller space} ${\cal T}(S)$ 
of $S$ is a smooth manifold diffeomorphic to 
$\mathbb{R}^{6g-6+2m}$. It can be represented as
the quotient of the space of all
complete hyperbolic metrics on $S$ of finite volume
under the action of the group of all diffeomorphisms
of $S$ which are isotopic to the identity.
The 
group of all isotopy classes of orientation preserving
diffeomorphisms of $S$ is
called the \emph{mapping class group} ${\cal M}(S)$ of $S$. 
It acts smoothly and
properly discontinuously on ${\cal T}(S)$ preserving
a complete Finsler metric, the so-called
\emph{Teichm\"uller metric}. The quotient
orbifold ${\rm Mod}(S)={\cal T}(S)/{\cal M}(S)$,
equipped with the projection of the Teichm\"uller
metric, is a complete noncompact geodesic metric space.

Even though the Teichm\"uller metric
is not non-positively curved in any reasonable
sense, it shares many properties with a
Riemannian metric of non-positive curvature.
For example, any two points in 
${\cal T}(S)$ can be connected by a unique
Teichm\"uller geodesic, and \emph{closed}
geodesics in the moduli space ${\rm Mod}(S)$
are in one-to-one correspondence
with the conjugacy classes of the so-called
\emph{pseudo-Anosov} elements of the 
mapping class group. 
However,
unlike in the case of negatively curved
Riemannian manifolds of finite volume, closed geodesics
may escape into the end of ${\rm Mod}(S)$.
Namely, we show.

\bigskip

{\bf Theorem :} {\it If $3g-3+m\geq 4$ then 
for every compact subset $K$ of ${\rm Mod}(S)$
there is a closed Teichm\"uller geodesic which
does not intersect $K$.
}

\bigskip

The organization of this note is as follows.
In Section 2 we summarize the properties
of the train track complex 
${\cal T\cal T}$ which are needed for our purpose.
In Section 3 we use train tracks and splitting
sequences to show Theorem A.

\section{The complex of train tracks}

In this section we summarize some results
and constructions from
\cite{PH92,H05} which will be used throughout the paper
(compare also \cite{M03}).

Let $S$ be an
oriented surface of
genus $g\geq 0$ with $m\geq 0$ punctures and where $3g-3+m\geq 2$.
A \emph{train track} on $S$ is an embedded
1-complex $\tau\subset S$ whose edges
(called \emph{branches}) are smooth arcs with
well-defined tangent vectors at the endpoints. At any vertex
(called a \emph{switch}) the incident edges are mutually tangent.
Through each switch there is a path of class $C^1$
which is embedded
in $\tau$ and contains the switch in its interior. In
particular, the branches which are incident
on a fixed switch are divided into
``incoming'' and ``outgoing'' branches according to their inward
pointing tangent at the switch. Each closed curve component of
$\tau$ has a unique bivalent switch, and all other switches are at
least trivalent.
The complementary regions of the
train track have negative Euler characteristic, which means
that they are different from discs with $0,1$ or
$2$ cusps at the boundary and different from
annuli and once-punctured discs
with no cusps at the boundary.
We always identify train
tracks which are isotopic.

A \emph{trainpath} on a train track $\tau$ is a $C^1$-immersion
$\rho:[m,n]\to \tau\subset S$ which maps each interval $[k,k+1]$
$(m\leq k\leq n-1)$ onto a branch of $\tau$. The integer $n-m$ is
then called the \emph{length} of $\rho$. We sometimes identify a
trainpath on $S$ with its image in $\tau$. Each complementary
region of $\tau$ is bounded by a finite number of trainpaths which
either are simple closed curves or terminate at the cusps of the
region.
A \emph{subtrack} of a train track $\tau$ is a subset $\sigma$ of
$\tau$ which itself is a train track. Thus every switch of
$\sigma$ is also a switch of $\tau$, and every branch of $\sigma$
is an embedded trainpath of $\tau$. We write $\sigma<\tau$ if
$\sigma$ is a subtrack of $\tau$.

A train track is called \emph{generic} if all switches are
at most trivalent.
The train track $\tau$ is called \emph{transversely recurrent} if
every branch $b$ of $\tau$ is intersected by an embedded simple
closed curve $c=c(b)\subset S$ which intersects $\tau$
transversely and is such that $S-\tau-c$ does not contain an
embedded \emph{bigon}, i.e. a disc with two corners at the
boundary.

A \emph{transverse measure} on a train track $\tau$ is a
nonnegative weight function $\mu$ on the branches of $\tau$
satisfying the \emph{switch condition}:
For every switch $s$ of $\tau$, the sum of the weights
over all incoming branches at $s$
is required to coincide with the sum of
the weights over all outgoing branches at $s$.
The train track is called
\emph{recurrent} if it admits a transverse measure which is
positive on every branch. We call such a transverse measure $\mu$
\emph{positive}, and we write $\mu>0$.
If $\mu$ is any transverse measure on a train track
$\tau$ then the subset of $\tau$ consisting of all
branches with positive $\mu$-weight is a recurrent
subtrack of $\tau$.
A train track $\tau$ is called \emph{birecurrent} if
$\tau$ is recurrent and transversely recurrent.

A \emph{geodesic lamination} for a complete
hyperbolic structure on $S$ of finite volume is
a \emph{compact} subset of $S$ which is foliated into simple
geodesics.
A geodesic lamination $\lambda$ is called \emph{minimal}
if each of its half-leaves is dense in $\lambda$. Thus a simple
closed geodesic is a minimal geodesic lamination. A minimal
geodesic lamination with more than one leaf has uncountably
many leaves and is called \emph{minimal arational}.
Every geodesic lamination $\lambda$ consists of a disjoint union of
finitely many minimal components and a finite number of isolated
leaves. Each of the isolated leaves of $\lambda$ either is an
isolated closed geodesic and hence a minimal component, or it
\emph{spirals} about one or two minimal components
\cite{CEG87,O96}.

A geodesic
lamination is \emph{finite} if it contains only finitely many
leaves, and this is the case if and only if each minimal component
is a closed geodesic. A geodesic lamination is \emph{maximal}
if its complementary regions are all ideal triangles
or once punctured monogons.
A geodesic lamination $\lambda$
is called \emph{complete} if $\lambda$ is maximal and
can be approximated in the \emph{Hausdorff topology} by
simple closed geodesics. The space ${\cal C\cal L}$
of all complete geodesic laminations equipped with
the Hausdorff topology is a compact metrizable space.
Every geodesic lamination $\lambda$
which is a disjoint union of finitely many minimal components
is a \emph{sublamination} of
a complete geodesic lamination, i.e. there
is a complete geodesic lamination which contains
$\lambda$ as a closed subset \cite{H05}.

A train track or a geodesic lamination $\sigma$ is
\emph{carried} by a transversely recurrent train track $\tau$ if
there is a map $F:S\to S$ of class $C^1$ which is isotopic to the
identity and maps $\sigma$ into $\tau$ in such
a way that the restriction of the differential
of $F$ to the tangent space of $\sigma$ vanishes
nowhere. Note that this makes sense since
a train track has a tangent line everywhere.
We call the restriction of $F$ to
$\sigma$ a \emph{carrying map} for $\sigma$. Write $\sigma\prec
\tau$ if the train track or the
geodesic lamination $\sigma$ is carried by the train track
$\tau$. If $\sigma$ is a train track which is
carried by $\tau$, then every 
geodesic lamination $\lambda$ which is carried
by $\sigma$ is also carried by $\tau$.
A train track
$\tau$ is called \emph{complete} if it is generic and
transversely recurrent and if
it carries
a complete geodesic lamination. Every complete train track
is recurrent. The set
of all complete geodesic laminations 
which are carried by a complete
train track $\tau$ is open and closed in ${\cal C\cal L}$.
In particular, the space ${\cal C\cal L}$
is totally disconnected \cite{H05}.

A half-branch $\hat b$ in a generic train track $\tau$ incident on
a switch $v$ of $\tau$ is called
\emph{large} if every trainpath containing $v$ in its interior
passes through $\hat b$. A half-branch which is not large
is called \emph{small}.
A branch
$b$ in a generic train track
$\tau$ is called
\emph{large} if each of its two half-branches is
large; in this case $b$ is necessarily incident on two distinct
switches, and it is large at both of them. A branch is called
\emph{small} if each of its two half-branches is small. A branch
is called \emph{mixed} if one of its half-branches is large and
the other half-branch is small (for all this, see \cite{PH92} p.118).

There are two simple ways to modify a complete train track $\tau$
to another complete train track. First, we can \emph{shift}
$\tau$ along a mixed branch to a train track $\tau^\prime$ as shown
in Figure A below. If $\tau$ is complete then the same is true for
$\tau^\prime$. Moreover, a train track or a
lamination is carried by $\tau$ if and
only if it is carried by $\tau^\prime$ (see \cite{PH92} p.119).
In particular, the shift $\tau^\prime$ of $\tau$ is
carried by $\tau$. Note that there is a natural
bijection of the set of branches of $\tau$ onto
the set of branches of $\tau^\prime$.

\begin{figure}[ht]
\includegraphics{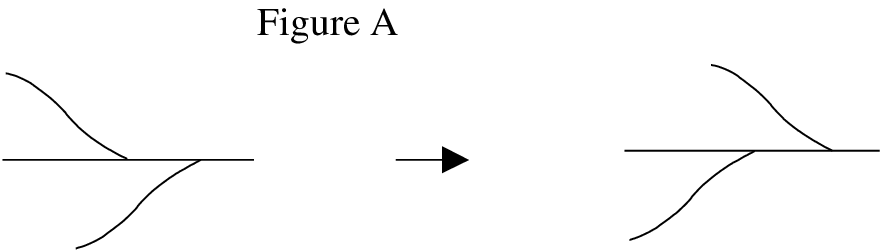}
\end{figure}

Second, if $e$ is a large branch of $\tau$ then we can perform a
right or left \emph{split} of $\tau$ at $e$ as shown in Figure B.
Note that a right split at $e$ is uniquely
determined by the orientation of $S$ and does not
depend on the orientation of $e$.
Using the labels in the figure, in the case of a right
split we call the branches $a$ and $c$ \emph{winners} of the
split, and the branches $b,d$ are \emph{losers} of the split. If
we perform a left split, then the branches $b,d$ are winners of
the split, and the branches $a,c$ are losers of the split.
The split $\tau^\prime$ of a train track $\tau$ is carried
by $\tau$, and there is a natural choice of a carrying map which
maps the switches of $\tau^\prime$ to the switches of $\tau$. The
image of a branch of $\tau^\prime$ is then a trainpath on $\tau$
whose length either equals one or two. We call
this carrying map the \emph{canonical carrying map}.
It induces
a natural bijection of the set of branches
of $\tau$ onto the set of branches of $\tau^\prime$ which
maps the branch $e$ to the \emph{diagonal} $e^\prime$ of the split.
The split of a maximal
transversely recurrent generic train track is maximal,
transversely recurrent and generic. If $\tau$ is complete and if
$\lambda\in {\cal C\cal L}$ is carried by $\tau$, then there is a
unique choice of a right or left split of $\tau$ at $e$ with the
property that the split track $\tau^\prime$ carries $\lambda$.
We call such a split a \emph{$\lambda$-split}.
The train track $\tau^\prime$ is
complete. In particular, a complete train track $\tau$ can always
be split at any large branch $e$ to a complete train track
$\tau^\prime$; however there may be a choice of a right or left
split at $e$ such that the resulting
train track is not recurrent any
more (compare p.120 in \cite{PH92}). The reverse of a split is called a
\emph{collapse}.

\begin{figure}[ht]
\includegraphics{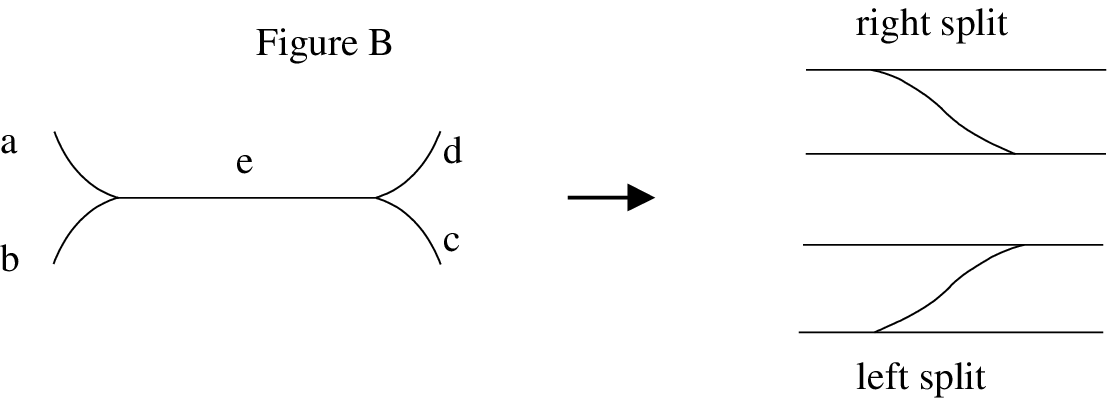}
\end{figure}

Denote by ${\cal T\cal T}$ the directed graph whose
vertices are the isotopy classes of complete train tracks on $S$
and whose edges are determined as follows. The train track
$\tau\in {\cal T\cal T}$
is connected to the train track $\tau^\prime$
by a directed edge if and only $\tau^\prime$ can be obtained
from $\tau$ by a single split.
The graph ${\cal T\cal T}$ is connected \cite{H05}.
As a consequence, if we identify each edge in ${\cal T\cal T}$
with the unit interval $[0,1]$ then this provides
${\cal T\cal T}$ with the structure of a connected locally finite metric
graph. Thus ${\cal T\cal T}$ is a locally compact complete geodesic
metric space. In the sequel we always assume that
${\cal T\cal T}$
is equipped with this metric without further comment.
The mapping class group ${\cal M}(S)$ of $S$ acts properly
and cocompactly on ${\cal T\cal T}$ as a group of
isometries. In particular, ${\cal T\cal T}$
is ${\cal M}(S)$-equivariantly quasi-isometric to
${\cal M}(S)$ equipped with any word metric \cite{H05}.

\section{Closed Teichm\"uller geodesics}

In this section we use train tracks and
splitting sequences to investigate closed
\emph{Teichm\"uller geodesics} in \emph{moduli space} 
${\rm Mod}(S)$ of our nonexceptional 
surface $S$ of finite type and show 
Theorem A from the introduction.
Such closed geodesics are the projections
of those Teichm\"uller geodesics
in the Teichm\"uller space
${\cal T}(S)$ which are invariant under
a pseudo-Anosov element of 
${\cal M}(S)$. Recall that
since the Euler characterstic of $S$ is negative
by assumption, ${\cal T}(S)$ can be identified with the
space of all  
marked complete hyperbolic metrics on $S$ of finite 
volume (with the usual identification of the elements
contained in an orbit for the action of the group
of diffeomorphisms of $S$ isotopic to the
identity).

A \emph{measured geodesic lamination} is 
a geodesic lamination $\lambda$ with a
translation invariant transverse measure $\nu$ such that
the $\nu$-weight of every compact
arc in $S$ with endpoints in $S-\lambda$ 
which intersects $\lambda$ nontrivially
and transversely is positive. We say
that $\lambda$ is the \emph{support} of the 
measured geodesic lamination. The space ${\cal M\cal L}$
of measured geodesic laminations 
equipped with the weak$^*$-topology admits a
natural continuous action of the multiplicative group
$(0,\infty)$. The quotient under this action is
the space ${\cal P\cal M\cal L}$ of 
\emph{projective measured geodesic laminations}
which is homeomorphic to the sphere $S^{6g-7+2m}$.
Every simple closed geodesic $c$ on $S$ defines
a measured geodesic lamination. The geometric
intersection number between simple closed
curves on $S$ extends to a continuous 
bilinear form $i$ on ${\cal M\cal L}$, the
\emph{intersection form}. We say that
a pair $(\lambda,\mu)\in {\cal M\cal L}\times
{\cal M\cal L}$ of measured geodesic laminations
\emph{jointly fills up} $S$ if 
for every measured geodesic
lamination $\eta\in {\cal M\cal L}$ we have
$i(\eta,\lambda)+i(\eta,\mu)>0$. This
is equivalent to saying that every complete 
simple (possibly infinite)
geodesic on $S$ intersects either the support
of $\lambda$ or the support of $\nu$ transversely.

The unit cotangent bundle of the Teichm\"uller space
with respect to the Teichm\"uller metric
can naturally be identified with 
the space ${\cal Q}^1(S)$ of unit area marked holomorphic
\emph{quadratic differentials} on our surface $S$.
Such a
quadratic differential $q$ is determined by a pair
$(\lambda^+,\lambda^-)$ of 
measured geodesic laminations on $S$
which jointly fill up $S$ and such that
$i(\lambda^+,\lambda^-)=1$ by our area normalization. 
The \emph{horizontal} measured
lamination $\lambda^+$ for $q$
corresponds to the equivalence class of the horizontal measured
foliation of $q$ (compare \cite{L83} and \cite{K92} for a discussion of
the relation between measured geodesic laminations and equivalence
classes of measured foliations on $S$). 
The quadratic differential $q$ determines a singular
euclidean metric on $S$ of unit area. For
every simple closed curve $c$ on $S$, the \emph{$q$-length} 
of $c$ is defined to be the infimum of the lengths
with respect to this metric of any curve which is
freely homotopic to $c$.

The bundle ${\cal Q}^1(S)$
admits a smooth action of the group $SL(2,\mathbb{R})$. The
restriction of this action to the one-parameter subgroup of
diagonal matrices of $SL(2,\mathbb{R})$ is the \emph{Teichm\"uller
geodesic flow} $\Phi^t$ on ${\cal Q}^1(S)$. If a quadratic
differential $q$ is given by the pair $(\lambda^+,\lambda^-)$ of
transverse measured geodesic laminations 
with $i(\lambda^+,\lambda^-)=1$ then $\Phi^tq$ is given
by the pair $(e^t\lambda^+,e^{-t}\lambda^-)$. The action of
$SL(2,\mathbb{R})$ commutes with the natural action of the mapping
class group ${\cal M}(S)$ on ${\cal Q}^1(S)$ and hence the flow
$\Phi^t$ projects to a flow on the quotient ${\cal Q}(S)={\cal
Q}^1(S)/{\cal M}(S)$, again denoted by $\Phi^t$. 
A flow line of this flow projects
to a Teichm\"uller geodesic in moduli space ${\cal T}(S)/
{\cal M}(S)={\rm Mod}(S)$. 

For every complete train track $\tau$, the
convex cone ${\cal V}(\tau)$ of all transverse
measures on $\tau$ 
can naturally be identified with the
set of all measured geodesic laminations
whose support is carried by $\tau$.
A \emph{tangential measure} on a complete
train track $\tau$ assigns to a branch $b$ of $\tau$
a weight $\mu(b)\geq 0$ such that for every
complementary triangle $T$ of $\tau$ with sides
$c_1,c_2,c_3$ we have $\mu(c_i)\leq \mu(c_{i+1})+
\mu(c_{i+2})$ (here indices are taken modulo 3).
Let ${\cal V}^*(\tau)$ be the convex cone of
all tangential measures on $\tau$.
By the results from Section 3.4 of \cite{PH92},
there is a one-to-one correspondence between
tangential measures on $\tau$ 
and measured geodesic laminations which
\emph{hit $\tau$ efficiently}, i.e.
measured geodesic laminations whose support $\lambda$
intersects $\tau$ transversely up to isotopy 
and is such that
$\tau\cup \lambda$ does not contain any embedded bigon.
In the sequel we often identify a measure
$\mu\in {\cal V}(\tau)$ (or $\nu\in {\cal V}^*(\tau)$)
with the measured geodesic lamination it defines.
With this identification, 
the bilinear pairing $<,>:{\cal V}(\tau)\times
{\cal V}^*(\tau)\to [0,\infty)$ defined by $<\mu,\nu>=\sum_b
\mu(b)\nu(b)$ is just the restriction of the
intersection form on ${\cal M\cal L}\times {\cal M\cal L}$
(see Section 3.4 of \cite{PH92}).

Denote by ${\cal B}(\tau)$ 
the set of all pairs $(\lambda,\nu)\in {\cal V}(\tau)\times
{\cal V}^*(\tau)$ where
$\lambda\in {\cal V}(\tau)$ is a transverse measure 
on $\tau$ of total weight
one and where $\nu\in {\cal V}^*(\tau)$ is 
a tangential measure  with
$<\lambda,\nu>=1$ and such that
$\lambda,\nu$ jointly fill up $S$. Every pair $(\lambda,\nu)\in 
{\cal B}(\tau)$ defines a quadratic
differential $q(\lambda,\nu)$ of area one.

If $\{\tau(i)\}$ is any splitting sequence, 
if $(\lambda,\nu)\in {\cal B}(\tau(0))$ 
and if $\lambda$ is carried
by $\tau(s)$ for some $s>0$ then for every
$i\in \{0,\dots,s\}$ there is a number $a(i)>0$
such that $(a(i)\lambda,a(i)^{-1}\nu)\in {\cal B}(\tau(i))$.
Denote by $[\lambda],[\nu]\in {\cal P\cal M\cal L}$ 
the projective
classes of $\lambda,\nu$ and for $i<s$ define
$\rho([\lambda],[\nu])(\tau(i))=a(i+1)/a(i)$.
We call $\rho$ a
\emph{roof function} for the pair $(\lambda,\nu)$. We have.

\begin{lemma}\label{lemma51}
Roof functions are uniformly
bounded.\end{lemma}
\begin{proof}
Let $\eta\in {\cal T\cal T}$ be
obtained from $\tau$ by a single split at a large
branch $e$ and let $\mu$ be a transverse measure
on $\eta$ of total weight one. Via the natural
carrying map $\eta\to \tau$, the measure $\mu$
defines a transverse measure $\mu_0$ on $\tau$.
By the definition of a roof function
we have to show that the
total weight of $\mu_0$ is bounded from above
by universal constant. However, this is immediate
from the fact that if $a,c$ are the losing branches
of the split connecting $\tau$ to $\eta$ and if
$e^\prime$ is the diagonal of the split 
then
$\mu_0(e)=\mu(e^\prime)+\mu(a)+\mu(c)$ and
$\mu_0(b)=\mu(b)$ for every branch $b\not=e$ of $\tau$ and
its corresponding branch of $\eta$.
This shows the lemma.
\end{proof}

For $\epsilon >0$ let ${\cal T}(S)_\epsilon$
be the $\epsilon$-thick part of Teichm\"uller
space consisting of all marked hyperbolic
metrics whose \emph{systole} (i.e. the length
of the shortest closed geodesic)
is at least $\epsilon$.
For sufficiently small $\epsilon$ the
set ${\cal T}(S)_\epsilon$ is
an ${\cal M}(S)$-invariant contractible
subset of ${\cal T}(S)$ on which the
mapping class group acts cocompactly. 
Let $P:{\cal Q}^1(S)\to {\cal T}(S)$
be the canonical projection and 
define
${\cal Q}^1(\epsilon)=\{q\in {\cal Q}^1(S)\mid 
Pq\in {\cal T}(S)_{\epsilon}\}$. The sets
${\cal Q}^1(\epsilon)$ are invariant
under the action of ${\cal M}(S)$. Their projections
${\cal Q}(\epsilon)={\cal Q}^1(\epsilon)/
{\cal M}(S)\subset {\cal Q}(S)$ 
to ${\cal Q}(S)$ are compact and satisfy
${\cal Q}(\epsilon)\subset {\cal Q}(\delta)$ for
$\epsilon>\delta$ and $\cup_{\epsilon >0}{\cal Q}(\epsilon)=
{\rm Mod}(S)$.

Call a finite splitting sequence $\{\tau(i)\}_{0\leq i\leq m}
\subset {\cal V}({\cal T\cal T})$ \emph{tight} if 
the natural carrying map $\tau(m)\to \tau(0)$
maps every branch $b$ of $\tau(m)$ \emph{onto} $\tau(0)$.
Note that every splitting sequence which
contains a tight subsequence is tight itself.
We have

\begin{lemma}\label{lemma52}
Let $\{\tau(i)\}_{0\leq i\leq m}\subset
{\cal V}({\cal T\cal T})$
by any finite splitting sequence such that
for some $\ell <m$, both sequences $\{\tau(i)\}_{0\leq i\leq\ell}$
and $\{\tau(i)\}_{\ell\leq i\leq m}$ are tight.
Then there is a number 
$\epsilon=\epsilon\{\tau(i)\}>0$ depending
on our sequence such that for all $(\lambda,\nu)\in 
{\cal B}(\tau(0))$ with the additional property that
$\lambda$ is carried by $\tau(m)$ the following holds.
\begin{enumerate}
\item $q(\lambda,\nu)\in {\cal Q}^1(\epsilon)$.
\item The minimal $\lambda$-weight of every branch of
$\tau(0)$ is at least $\epsilon$.
\end{enumerate}
\end{lemma}
\begin{proof}
Let $\{\tau(i)\}_{0\leq i\leq m}
\subset {\cal V}({\cal T\cal T})$ 
be a finite splitting sequence such that for
some $\ell >0$ the sequences $\{\tau(i)\}_{0\leq i\leq \ell}$
and $\{\tau(i)\}_{\ell\leq i\leq m}$ are tight. 
Extending the distance function on the \emph{curve
graph} of $S$, define a distance $d$ on the
space ${\cal U}$ of minimal geodesic laminations
which do \emph{not} fill up $S$
(i.e. which admit a complementary component which
is neither a disc nor a once punctured disc)
by requiring that
the distance between $\lambda\not=\mu$ equals one if
and only if $\lambda,\mu$ can be realized disjointly.
By the choice of ${\cal U}$, the distance
between any two elements $\alpha,\beta\in {\cal U}$
is finite. 
Let $\lambda\in {\cal U}$ be a minimal geodesic lamination
which is carried by $\tau(m)$ and let
$\beta\in {\cal U}$ be such that
$d(\beta,\lambda)=1$; we claim that $\beta$ is carried
by $\tau(\ell)$. Namely, since
our splitting sequence $\{\tau(i)\}_{\ell\leq i\leq m}$
is tight by assumption, the lamination
$\lambda$ fills $\tau(\ell)$, i.e. any
transverse measure for $\lambda$ defines
a \emph{positive} transverse measure on $\tau(\ell)$.
Assume first that $\beta$ is a simple closed curve
on $S$ with $d(\lambda,\beta)=1$.
Since $\lambda$ is minimal there is a sequence
$(\alpha_i)_i$ of simple
closed curves which converge to $\lambda$ in the
Hausdorff topology \cite{CEG87}. Now $\lambda$ is disjoint from
$\beta$ and therefore we can choose 
the sequence $(\alpha_i)_i$ in
such a way that each of the curves
$\alpha_i$ is disjoint from $\beta$.
For sufficiently large $i$, the curve
$\alpha_i$ is carried by $\tau(m)$ (see Lemma 2.4 of \cite{H05})
and fills $\tau(\ell)$, i.e. a
carrying map $\alpha_i\to \tau(\ell)$
is surjective. By Lemma 4.4 of
\cite{MM99}, since $\alpha_i,\beta$ are
disjoint the curve $\beta$ is carried by
$\tau(\ell)$. Now the space of geodesic
laminations carried by $\tau(\ell)$ is closed
with respect to the Hausdorff topology and therefore
$\tau(\ell)$ carries every geodesic lamination 
$\beta\in {\cal U}$ with
$d(\lambda,\beta)=1$.

Let again $\lambda$ be a minimal
geodesic lamination which is carried by $\tau(m)$
and let $\nu$ be a minimal geodesic lamination which
hits $\tau(0)$ efficiently. We claim that $\lambda$ and
$\nu$ jointly fill up $S$. Namely, otherwise
there is a simple closed curve $c$ which is disjoint
from both $\lambda,\nu$. But this means that
the distance between $\lambda$ and $c$ is at most one
and hence $c$ is carried by $\tau(\ell)$ by our
above consideration. The same argument, applied
to $c$ and $\nu$ and the
splitting sequence $\{\tau(i)\}_{0\leq i\leq \ell}$, 
shows that $\nu$ is carried
by $\tau(0)$ which contradicts our assumption
that $\nu$ hits $\tau(0)$ efficiently.

Let $\mu$ be any transverse measure
for $\tau(m)$ and let $m_0=\max\{\mu(b)\mid b$
is a branch of $\tau(m)\}$.
Let $\mu_0$ be the transverse measure
on $\tau(0)$ induced from $\mu$ 
by the carrying map
$\tau(m)\to \tau(0)$. Since the 
splitting sequence $\{\tau(i)\}_{0\leq i\leq m}$ is tight, 
the $\mu_0$-weight of a branch $b$ of $\tau(0)$
is not smaller than $m_0$, and it is not
bigger than $c m_0$ where $c>0$
only depends on our splitting sequence
$\{\tau(i)\}_{0\leq i\leq m}$
but not on $\mu$.
As a consequence, there is a universal constant
$s>0$ not depending on $\mu$ with the property that if 
the measure $\mu$ on $\tau(m)$ is
normalized in such a way that
the total weight of $\mu_0$
on $\tau(0)$ equals one, then the minimum
of the weights that the measure $\mu_0$ disposes
on the branches of $\tau(0)$ is not smaller
than $s$. This shows the
second statement in our lemma. Now for any tangential measure
$\nu$ on $\tau(0)$ we have
$i(\mu,\nu)=\sum_b\mu_0(b)\nu(b)$ (where as before we
identify the measures $\mu,\nu$ with the 
measured geodesic laminations they define)
and therefore if $i(\mu,\nu)=1$ then the
\emph{maximal} weight
disposed on a branch of $\tau(0)$ by the
tangential measure $\nu$ 
is not bigger than $1/s$. Since ${\cal V}(\tau(m))$
and ${\cal V}^*(\tau(0))$ are closed subsets
of ${\cal M\cal L}$ and the intersection
form is continuous, this implies that 
the set ${\cal B}_0$ of pairs 
$(\mu,\nu)\in {\cal V}(\tau(m))\times 
{\cal V}^*(\tau(0))\subset {\cal V}(\tau(0))
\times {\cal V}^*(\tau(0))$ such that
$\mu$ defines a
transverse measure of total weight one on $\tau(0)$ and
such that $i(\mu,\nu)=1$ is a compact
subset of ${\cal V}(\tau(0))\times {\cal V}^*(\tau(0))$.
By our above consideration, every pair 
$(\mu,\nu)\in {\cal B}_0$ jointly fills up $S$
and hence ${\cal B}_0$ is a compact subset
of ${\cal B}(\tau(0))$.

Since the assignment which associates
to a quadratic differential
its horizontal and vertical measured foliation
is injective, continuous and open with
respect to the smooth topology
on ${\cal Q}^1(S)$ and the weak topology
on ${\cal M\cal L}\times{\cal M\cal L}$
the set
${\cal C}=\{q(\mu,\nu)\mid (\mu,\nu)\in {\cal B}_0\}$
is compact. Therefore there is some
$\epsilon >0$ such that ${\cal C}\subset
{\cal Q}^1(\epsilon)$ which shows the lemma.
\end{proof}

For the proof of the theorem
from the introduction, we use
biinfinite splitting 
sequences $\{\tau(i)\}_i\subset {\cal V}({\cal T\cal T})$ 
to construct
closed orbits of the Teichm\"uller flow.
For this recall that there are only finitely
many orbits of complete train tracks under
the action of ${\cal M}(S)$. 
Call a biinfinite splitting sequence
$\{\tau(i)\}$ \emph{periodic} if there is a number
$m>0$ and some $\phi\in {\cal M}(S)$ 
such that for every $k\in \mathbb{Z}$ we have
$\{\tau(i)\}_{km\leq i\leq (k+1)m}=
\phi^k\{\tau(i)\}_{0\leq i\leq m}$.
We then call the mapping class $\phi$ a
\emph{period map} of the sequence with \emph{period}
$m$. A periodic splitting
sequence can easily be constructed 
by choosing any infinite splitting sequence
$\{\eta(i)\}$ and some $i<j$ such that
$\eta(j)=\phi\eta(i)$ for some
$\phi\in {\cal M}(S)$; then 
$\tau(k(j-i)+s)=\phi^k\eta(i+s)$ $(s<j-i,k\in \mathbb{Z})$
defines a periodic splitting sequence with
period map $\phi$.
Note that a period map of any periodic
splitting sequence 
is an element of ${\cal M}(S)$ of infinite order.
We say that a measured geodesic lamination
$\lambda$ \emph{fills up} $S$ if for every
simple closed curve $c$ on $S$ we have $i(c,\lambda)>0$.
The support of such a measured geodesic
lamination is necessarily connected. 
If $\lambda\in {\cal M\cal L}$ is such 
that its projectivization $[\lambda]\in {\cal P\cal M\cal L}$
is a fixed point of a pseudo-Anosov mapping class
then $\lambda$ fills up $S$.
Recall that for every complete train track $\tau$ 
we view ${\cal V}(\tau)$ as the
space of all measured geodesic laminations which
are carried
by $\tau$. We have

\begin{lemma}\label{lemma53}
Let $\{\tau(i)\}
\subset {\cal V}({\cal T\cal T})$ be
a biinfinite periodic splitting sequence
with period map $\phi\in {\cal M}(S)$. Then
$\phi$ is pseudo-Anosov if and only if
there is some 
$\lambda\in \cap_i {\cal V}(\tau(i))$
which fills up $S$.
\end{lemma}
\begin{proof}
Let $\{\tau(i)\}_i$ be a biinfinite
periodic splitting sequence with 
period map $\phi\in {\cal M}(S)$. 
Then there is a number $\ell>0$ such that
$\phi{\cal V}(\tau(0))={\cal V}(\tau(\ell))\subset
{\cal V}(\tau(0))$ and $\phi^{-1}{\cal V}^*(\tau(0))=
{\cal V}^*(\tau(-\ell))\subset
{\cal V}^*(\tau(0))$.
Assume first that
$\phi$ is pseudo-Anosov. Since $\phi$ acts on
the Thurston boundary ${\cal P\cal M\cal L}$ 
of Teichm\"uller space
with north-south-dynamics and since the projectivizations
$P,Q$ of the cones ${\cal V}(\tau(0)),{\cal V}^*(\tau(0))$ 
are closed disjoint subsets of ${\cal P\cal M\cal L}$ 
with non-empty
interior, the attracting fixed point $\lambda^+$ 
for the action of $\phi$ on ${\cal P\cal M\cal L}$
is contained in $\cap_i\phi^i(P)$. But 
$\lambda^+$ is 
the projectivization of a measured geodesic lamination
$\lambda$ which fills up $S$ 
and hence $\lambda\in \cap_i{\cal V}(\tau(i))$ is
a lamination as required in the lemma.

Now assume that there is some
$\lambda\in \cap_i{\cal V}(\tau(i))=
\cap_i\phi^i{\cal V}(\tau(0))=A$ which
fills up $S$. The set $A$ is a closed
non-empty $\phi$-invariant subset of ${\cal M\cal L}$.
By Theorem 8.5.1 of \cite{M03}, it 
can be represented 
as a join ${\cal V}_1 * {\cal V}_2$ where
${\cal V}_1\not=\emptyset$ 
is the space of measured geodesic laminations
whose support is contained in a fixed
geodesic lamination $\zeta$ on $S$ and where
${\cal V}_2$ is a space of measured geodesic laminations
whose support is carried by a train track on $S$ contained
in the complement of $\zeta$. A measured geodesic
lamination whose support contains more than one
connected component does not fill up $S$.
But $\lambda\in {\cal V}_1*{\cal V}_2$ fills up
$S$ by assumption 
and therefore the set ${\cal V}_2$ 
is empty and $A$ is the
space of measured geodesic laminations supported in the
single geodesic lamination $\zeta$. Moreover,
$\zeta$ is minimal and connected, 
with complementary components which
are topological discs or once punctured
topological discs. Every
mapping class which preserves the set $A$ also
preserves the geodesic lamination
$\zeta$ and hence it is pseudo-Anosov
by Thurston's classification of elements of 
${\cal M}(S)$.
But $A$ is invariant under the period map 
$\phi$ and hence
$\phi$ is pseudo-Anosov as claimed.
\end{proof}

Let $p>0$ be the number of branches of 
a complete train track on $S$. 
Let $\tau\in {\cal V}({\cal T\cal T})$ and 
number
the branches of $\tau$ with numbers
$1,\dots,p$. If 
$\tau^\prime$ can be obtained from $\tau$ by a
single split at a large branch $e$, then
the numbering of the branches of 
$\tau$ naturally induces a numbering
of the branches of $\tau^\prime$
(compare the discussion in Section 5 of \cite{H05}).
In other words, every splitting sequence
$\{\tau(i)\}$ together with a numbering
of the branches of $\tau(0)$ determines
uniquely a \emph{numbered splitting sequence},
i.e. a splitting sequence together with a
consistent numbering of the branches of
the train tracks in the sequence which
is determined by 
$\{\tau(i)\}$ up to a permutation of the
numbering of the branches of $\tau(0)$.
Let $e_1,\dots,e_p$ be the standard basis
of $\mathbb{R}^p$. A numbering of the
branches of $\tau(0)$ then defines
an embedding of ${\cal V}(\tau(0))$ onto
a closed convex cone in $\mathbb{R}^p$ determined
by the switch conditions by associating 
to a measure $\mu\in {\cal V}(\tau(0))$ 
the vector $\sum_i\mu(i)e_i\in \mathbb{R}^q$ where we
identify a branch of $\tau(0)$ with its number.
If $\tau^\prime$ is obtained
from $\tau$ by a split at the large branch $e$
with number $k$
and if $i,j$ are the losing branches of the split
connecting $\tau$ to $\tau^\prime$
then the transformation ${\cal V}(\tau^\prime)\to
{\cal V}(\tau)$ induced by the canonical carrying
map $\tau^\prime\to \tau$ is just the 
restriction to ${\cal V}(\tau^\prime)$ of the 
unique linear map
$A:\mathbb{R}^p\to \mathbb{R}^p$ 
which satisfies $A(e_i)=e_i+e_k,
A(e_j)=e_j+e_k$ and $A(e_s)=e_s$ for $s\not=i,j$.
A change of the numbering of the branches
of $\tau(0)$ results in replacing
$A$ by its conjugate under a \emph{permutation map}, i.e.
a linear
isometry of $\mathbb{R}^q$ corresponding to
a permutation of the standard basis vectors.
Similarly, for every finite splitting sequence
$\{\tau(i)\}_{0\leq i\leq m}$ the canonical map
${\cal V}(\tau(m))\to {\cal V}(\tau(0))$ 
corresponds to a $(p,p)$-matrix $A(\tau(m))$ 
with non-negative entries
which is uniquely determined by the sequence up to
conjugation with a permutation matrix.

By the Perron Frobenius theorem,
a $(p,p)$-matrix $A$ with positive entries
admits up to a multiple a
unique eigenvector with positive
entries. 
The corresponding eigenvalue $\alpha$
is positive, and its absolute value is
bigger than the absolute value of 
any other eigenvalue of $A$. 
Moreover, the generalized eigenspace for
$\alpha$ is one-dimensional.
We call an eigenvector with positive entries
for the eigenvalue
$\alpha$ a \emph{Perron Frobenius
eigenvector}. We have.

\begin{lemma}\label{lemma54}
Let $\{\tau(i)\}$ be 
a periodic splitting sequence with period map $\phi$.
Assume that $k>0$ is such that $\phi(\tau(0))=\tau(k)$
and that $\{\tau(i)\}_{0\leq i\leq k}$ is tight; 
then $\phi$ is pseudo-Anosov
and there is some $\ell \geq 1$ such that
its attracting fixed point is the projectivization
of a Perron Frobenius eigenvector of the matrix $A(\tau(k\ell))$.
\end{lemma}
\begin{proof}
Let $\{\tau(i)\}_{0\leq i\leq k}$
be a tight splitting sequence and assume that there is some
$\phi\in {\cal M}(S)$ such that $\phi(\tau(0))=\tau(k)$. We
obtain a biinfinite splitting sequence $\{\tau(i)\}$ 
by defining
$\tau(mk+\ell)=\phi^m(\tau(\ell))$ for $m\in \mathbb{Z}$ and 
$\ell <k$. By our above discussion, with
respect to some numbering of the branches of 
$\tau(0)$ the self-map ${\cal V}(\tau(0))\to
{\cal V}(\tau(0))$ which is the composition of 
the map ${\cal V}(\tau(0))\to {\cal V}(\tau(k))$ 
induced by $\phi$ and the transformation 
${\cal V}(\tau(k))\to {\cal V}(\tau(0))$  
induced by 
the canonical carrying map $\tau(k)\to \tau(0)$ 
can be represented by the restriction of a linear
map $\mathbb{R}^p\to \mathbb{R}^p$ given
with respect to the standard basis by a matrix
$C$. Since $\{\tau(i)\}_{0\leq i\leq k}$ is 
tight by assumption, the entries of the matrix
$C$ are positive and therefore by the
Perron Frobenius theorem, our map uniformly
contracts the projectivization of the closed cone
of vectors with nonnegative entries into its
interior. As a consequence,
the intersection $\cap_i{\cal V}(\tau(i))$
consists of a single ray spanned by the positive measure
on $\tau(0)$ which corresponds to a Perron Frobenius
eigenvector of $C$ 
(compare with the beautiful
argument in \cite{K85}). This ray defines a projective
measured geodesic lamination $\lambda^+
\in {\cal P\cal M\cal L}$ which is \emph{uniquely ergodic}
and fills up $S$.
Similarly, the intersection $\cap{\cal V}^*(\tau(i))$ consists
of a single ray which defines a projective
measured geodesic lamination $\lambda^-\not=\lambda^+$ 
which fills up $S$. By Lemma 5.3,
the mapping class
$\phi$ is pseudo-Anosov and its attracting fixed point
equals $\lambda^+$, its repelling fixed point
equals $\lambda^-$.

Recall that a numbering of the branches of $\tau(0)$
induces 
a numbering of the branches of $\tau(k)$. 
There is a second numbering of the branches of
$\tau(k)$ induced from the numbering of the branches
of $\tau(0)$ via the map $\phi:\tau(0)\to \tau(k)$.
These two numberings differ by 
a permutation $\sigma$ of $\{1,\dots,p\}$.
Let $\ell \geq 1$ be such that
$\sigma^\ell=1$; then $\phi^\ell\tau(0)=\tau(k\ell)$ as numbered train
tracks. It is then immediate from our above
discussion that the attracting 
fixed point of $\phi$ corresponds to the projectivization
of a Perron Frobenius eigenvector for the matrix
$A(\tau(k\ell))$ which shows the
lemma.
\end{proof}

In a Riemannian manifold $M$ of bounded negative
curvature and finite volume
with fundamental group $\Gamma$, every
closed geodesic intersects a fixed compact
subset of $M$. Moreover, we can construct
an infinite sequence $(\gamma_i)$ of
pairwise distinct  
closed geodesics which intersect every
compact subset $K$ of $M$ in arcs of uniformly
bounded length as follows. Choose a hyperbolic
element $\phi$ in $\Gamma$ and a parabolic
element $\psi$; then up to replacing $\psi$ by
a conjugate, the elements $\phi\psi^k$ are
hyperbolic for all $k>0$ and they define
an infinite family of pairwise distinct closed
geodesics with the required properties.
In the next lemma, we point out that
a similar construction for a pseudo-Anosov element
of ${\cal M}(S)$ and a Dehn twist about a
suitably chosen simple closed curve
(viewed as a parabolic element in ${\cal M}(S)$)
yields sequences of closed geodesics in
${\rm Mod}(S)$ with similar properties.

\begin{lemma}\label{lemma55}
There is a 
sequence $(\gamma_i)$ of pairwise distinct closed
Teichm\"uller geodesics 
in moduli space ${\rm Mod}(S)$ with the following
properties.
\begin{enumerate}
\item There is a fixed compact subset $K_0$
of ${\rm Mod}(S)$ which is intersected by
$\gamma_i$ for every $i$.
\item The geodesics $\gamma_i$ 
intersect every 
compact subset $K\subset {\rm Mod}(S)$
in arcs of uniformly bounded length.
\end{enumerate}
\end{lemma}
\begin{proof}
Recall from Section 3 the definition
of a train track $\sigma$ in special standard form
for a framing $F$ of our surface $S$ with 
pants decomposition $P$.
We first claim that there
is a train track $\tau$ which is shift equivalent
to a such a train track in special standard
form and there is
a tight splitting sequence $\{\tau(i)\}_{0\leq i\leq k}$
issuing from $\tau(0)=\tau$ such that $\tau(k)=\phi
\tau(0)$ for some $\phi\in {\cal M}(S)$.

For this let $\eta$ by 
any train track which is
shift equivalent to a train track in special
standard form for a framing $F$ of $S$ 
with pants decomposition $P$. Let $\lambda_0$
be the unique special geodesic lamination
for $P$ which is carried by $\eta$ and 
let $\{\zeta(i)\}_{0\leq i\leq m}$ be
any tight splitting sequence issuing from
$\zeta(0)=\eta$. Since every orbit for the
action of the mapping class group on the
space ${\cal C\cal L}$ of complete geodesic
laminations on $S$ is dense and since
the set of all complete geodesic
laminations on $S$ which are carried
by $\zeta(m)$ is open in ${\cal C\cal L}$, 
there is an element $\tilde \phi_0\in {\cal M\cal S}$ such 
that the special geodesic 
lamination $\tilde \phi_0\lambda_0$ for the pants
decomposition $\phi_0P$ is carried by $\zeta(m)$.
By the considerations in Section 4, the train 
track $\zeta(m)$ is splittable to 
a train track $\eta_1$ which is shift equivalent to 
the train track 
$\phi_0\eta$ for some $\phi_0\in {\cal M}(S)$. 
Note that every splitting
sequence connecting $\eta$ to $\eta_1$
is tight and therefore if $\eta_1=\phi_0\eta$ 
then such a splitting sequence has the desired 
properties. Otherwise we apply this
construction to $\eta_1$ and find a tight splitting
sequence connecting $\eta_1$ to a train track
$\eta_2$ which is shift equivalent to $\phi_1\eta$
for some $\phi_1\in {\cal M}(S)$. Since there
are only finitely many train track in the
shift equivalence class of $\eta$, after
finitely many steps we find some $0\leq i<j$ 
such that $\eta_j=(\phi_j\circ \phi_i^{-1})\eta_i$
and hence $\eta_i$ is splittable to 
$(\phi_j\circ \phi_i^{-1})\eta_i$ with
a tight splitting sequence.

Now let $\{\tau(i)\}_{0\leq i\leq k}$ be a tight
splitting sequence
as above with $\tau(k)=\phi(\tau(0))$ for some
$\phi\in {\cal M}(S)$ and such 
that $\tau(0)$ is shift equivalent
to a train track in special standard form
for a framing $F$ with pants decomposition $P$. By Lemma 5.4,
$\phi$ is pseudo-Anosov. For $k\leq i\leq 2k$ define
$\tau(i)=\phi(\tau(i-k))$.
By the definition of a train track in special standard
form, every component $\alpha$
of $P$ defines an embedded simple closed
trainpath on $\tau(0)$, and there
is a number $s\geq 2$ (depending
on $\alpha$) and a
Dehn twist $\psi$ about $\alpha$ with the
property that $\psi\tau(0)$ is splittable
to $\tau(0)$ with a
splitting sequence of length $s$. Each split
of this sequence is a split at a large
branch contained in $\alpha$. In other words,
for every $u>0$ the train track
$\psi^u\tau(0)$ is splittable to 
$\phi^2(\tau(0))=\tau(2k)$ with a tight splitting sequence
$\{\tau(i)\}_{-su\leq i\leq 2k}$ 
of length $2k+su$. By Lemma 5.4, the mapping class
$\zeta_u=\phi^2\psi^{-u}$ is pseudo-Anosov. Its
axis is the Teichm\"uller geodesic which is
defined by a quadratic differential
whose horizontal foliation
is contained in the ray $\cap_i\zeta_u^i
{\cal V}(\tau(0))$ and whose vertical
foliation is contained in the ray
$\cap_i\zeta_u^{-i}{\cal V}^*(\tau(0))$.

For $u\geq 1$ let $(\lambda_u,\nu_u)\in {\cal B}(\tau(0))$ 
be such that the quadratic
differential $q_u=q(\lambda_u,\nu_u)$ 
is a cotangent vector of the axis of 
$\zeta_u$. For $i\in [-su,2k]$
let $a_u(i)>0$ be such that $(a_u(i)\lambda_u,
a_u(i)^{-1}\nu_u)\in {\cal B}(\tau(i))$, i.e. that
the total weight
of the transverse measure on $\tau(i)$ which corresponds to
the measured geodesic lamination 
$a_u(i)\lambda_u$ equals one. 
Then the function $i\to a_u(i)$ is increasing, and
$\log a_u(2k)-\log a_u(-su)$ is the length of the
closed Teichm\"uller geodesic $\gamma_u$ 
in moduli space which is the projection of 
the axis of $\zeta_u$.

By Lemma \ref{lemma51}, the numbers 
$\log a_u(2k)$ are bounded
independent of $u$, and 
by Lemma \ref{lemma52} there is a number $\epsilon_0>0$
such that $q_u\in {\cal Q}^1(\epsilon_0)$ 
for every $u>0$. The projection 
${\cal Q}(\epsilon_0)$ of ${\cal Q}^1(\epsilon_0)$
to ${\cal Q}(S)$ is compact
and therefore the geodesics $\gamma_u$ have property 1)
stated in the lemma. Moreover,
since $\{\tau(i)\}_{0\leq i\leq k}$
is tight by assumption, by Lemma 5.2 there is a number 
$\delta>0$ such that for \emph{every} 
$u>0$ the $\lambda_u$-weight of \emph{every} branch
$b$ of $\tau(0)$ is at least $\delta$; then
the maximal weight that the tangential measure
$\nu_u$ disposes on a branch of 
$\tau(0)$ is at most $1/\delta$.

To show that the geodesics $\gamma_u$ also
satisfy property 2) stated in the lemma,
we show that $a_u(-su)\to 0$ $(u\to \infty)$ and that
moreover for every $\epsilon >0$ there is a number
$\beta(\epsilon)>0$ which is independent of $u$
and such that for all 
$t\in [-\log a_u(-su)+\beta(\epsilon),-\beta(\epsilon)]$
the $\Phi^tq_u$-length of the curve $\alpha$
is not bigger than $\epsilon$. To show
that such a constant exists,
denote by $b_1,\dots,b_r$ the 
branches of $\tau(0)$ which are incident on
a switch in the embedded trainpath $\alpha$ 
but which are not contained
in $\alpha$. Note that $r\geq 2$, i.e. there are
at least two such branches.
By the considerations in the proof
of Lemma 2.5 of \cite{H04} we have
$i(\mu,\alpha)=\frac{1}{2}\sum_i
\mu(b_i)$ for every $\mu\in {\cal V}(\tau(0))$.
In particular, the intersection number
$i(\lambda_u,\alpha)$ is bounded from
above independent of $u$.
Choose a numbering of the branches of $\tau(0)$ so
that $b_i$ is the branch with number $i$ and
such that the branches with numbers $r+1,\dots,r+\ell$
are precisely the branches contained in
the embedded trainpath $\alpha$.
Then the linear self-map of $\mathbb{R}^p$ which 
defines the natural map ${\cal V}(\tau(0))\to
{\cal V}(\tau(-s))$ with respect to this numbering
maps for each $i\leq r$ the 
$i$-th standard basis
vector $e_i$ of $\mathbb{R}^p$
to $e_i+e_{r+1}+\dots +e_{r+\ell}$, and
for $i>r$ it maps $e_i$ to itself.
For $\delta>0$ as above,
the sum of the $\lambda_u$-weights of the branches
$b_i$ $(i\leq r)$ 
is bounded from below by $2\delta$ and
hence for every $m\leq u$
we have $a_u(-sm)\leq 1/(1+2m\delta)$ 
independent of $u$. In particular, $a_u(-su)$ tends
to zero as $u\to \infty$.
However, $i(c\lambda_u,\alpha)=
c i(\lambda_u,\alpha)$ for all $c>0$ and therefore
for every $\epsilon >0$ there is a number
$\beta_1(\epsilon)>0$ not depending on $u$ such that
for $t\in [\log a_u(-su), -\beta_1(\epsilon)]$ the intersection
with $\alpha$ 
of the horizontal measured geodesic
lamination of $\Phi^t q_u$ 
does not exceed $\epsilon/4$.

We observed above that the maximal weight
that the tangential measure $\nu_u$ disposes
on a branch of $\tau(0)$ is bounded from
above by $1/\delta$. Since the natural transformation
${\cal V}^*(\tau(0))\to {\cal V}^*(\tau(2k))$
can be represented by a fixed linear map
(with the interpretation discussed above) and the
functions $a_u$ are non-decreasing,
there is then a number
$\chi>0$ with the property that
for \emph{every} $u>0$ the maximal weight that
the measure $a_u(2k)^{-1}\nu_u$ disposes
on a branch of $\tau(2k)$ is bounded from
above by $\chi$. However, the tangential measure
$a_u(-su)^{-1}\nu_u$ on $\tau(-su)$ 
is the image of the tangential 
measure $a_u(2k)^{-1}\nu_u$ 
on $\tau(2k)$
under the transformation $\zeta_u^{-1}$ 
which identifies $\tau(2k)$ with $\tau(-su)$ 
and therefore the maximal weight that
$a_u(-su)^{-1}\nu_u$ disposes on
a branch of $\tau(-su)$ is bounded from above
by $\chi$. Since $\alpha$ is 
an embedded trainpath in 
$\tau(-su)$ this just means that
the intersection with $\alpha$ of the
vertical geodesic lamination of 
$\Phi^{\log a_u(-su)}q_u$ is bounded from
above by a universal constant not depending
on $u$. In other words, there is a constant
$\beta_2(\epsilon)>0$ not depending on $u$
such that for every $t\in 
[-\log a_u(-su)+\beta_2(\epsilon),0]$,
the intersection with $\alpha$ of the
vertical measured geodesic lamination of $\Phi^tq_u$ does
not exceed $\epsilon/4$. The
$\Phi^tq_u$-length of $\alpha$ is bounded from
above by twice the sum of the intersection
of $\alpha$ 
with the horizontal and the vertical
geodesic lamination of $\Phi^tq_u$ 
(compare e.g. \cite{R05}) and hence
for every sufficiently large $u$ and for every
$t\in [-\log a_u(-su)+\beta_2(\epsilon),
-\beta_1(\epsilon)]$ the $\Phi^tq_u$-length
of $\alpha$ does not exceed $\epsilon$.
Now for every $\delta >0$ there is a number
$\rho(\delta)>0$ such that the set
${\cal Q}^1(\delta)$ is contained in the
set of all quadratic differentials $q$ with the
property that the $q$-length of \emph{every}
essential simple closed
curve on $S$ is at least $\rho(\delta)$.
This then yields that the geodesics $\gamma_u$
intersect a fixed compact subset $K$ of 
${\rm Mod}(S)$ 
in arcs of uniformly bounded length, independent of $u$.
The lemma follows.
\end{proof}

We are now ready to show the
Theorem from the introduction.

\begin{proposition}\label{prop56}
If $3g-3+m\geq 4$ then 
for every $\epsilon >0$
there is a closed orbit of the Teichm\"uller
geodesic flow which does not intersect
${\cal Q}(\epsilon)$.\end{proposition}
\begin{proof} Assume that $3g-3+m\geq 4$ (which
excludes a sphere with at most 6
punctures, a torus with at most 
3 punctures and a closed surface of genus 2) and 
let $P$ be a pants decomposition of $S$.
We require that if the genus of $S$ is at least 2 then 
$P$ contains two non-separating
simple closed curves $\gamma_1,\gamma_2$ such that
the surface $S_0$ obtained by cutting $S$ open 
along $\gamma_1,\gamma_2$ is connected and that
moreover the bordered surface $S-(P-\gamma_1-\gamma_2)$
consists of $2g-6+m$ pairs of pants and 
2 forth punctured spheres containing
$\gamma_1,\gamma_2$ as essential curves in their interior.
If the genus of $S$ is at most one then $S$ contains
at least 4 punctures and 
we choose $P$ in such a way that
it contains 
two separating simple closed curves $\gamma_1,\gamma_2$
with the property that
the surface obtained by cutting $S$ open along
$\gamma_1,\gamma_2$ consists of three connected
components, where two of these components
are pairs of pants. In other words, $S-(\gamma_1 \cup
\gamma_2)$ contains a unique connected component
$S_0$ of Euler characteristic at most $-2$, and
if we replace the boundary circles of this component
by cusps then the resulting surface is
nonexceptional.

Let $\tau$ be a train track
in standard form for a framing with
pants decomposition $P$ and with only twist connectors.
Then the curves
$\gamma_1,\gamma_2$ are carried by $\tau$ and 
define embedded simple
closed trainpaths on $\tau$. We require that
a branch of $\tau$ which is contained in the closure
of $S_0$, which is incident on a switch contained
in $\gamma_1\cup \gamma_2$ and which is 
\emph{not} contained in $\gamma_1\cup\gamma_2$ 
is a small branch. By our choice of $P$ and 
the discussion on p.147-48 of \cite{PH92},
such a complete train track exists.
After removing
all branches of $\tau$ which are incident
on a switch contained in $\gamma_1\cup\gamma_2$ 
or which are contained in $S-S_0$ we obtain
a train track $\sigma_0$ on 
the subsurface $S_0$. If we identify $S_0$ with 
the surface of finite type obtained by replacing
the boundary circles by cusps then $\sigma_0$ is
a complete train track on $S_0$ which is 
in standard form for a framing with pants
decomposition $P\cap S_0$. 
Similarly, for $i=1,2$ let 
$S_i\supset S_0$ be the 
connected component of $S-\gamma_i$ which 
contains $S_0$ as a subsurface.
The train track $\sigma_i$ on $S_i$ 
which we obtain from $\tau$ by removing all
branches which are incident on a switch in $\gamma_i$ 
or which are not contained in $S_i$ 
is complete (in the interpretation 
as above). Note that we have
$\tau=\sigma_1\cup \sigma_2$, i.e. every branch of
$\tau$ is either a branch of $\sigma_1$ or a branch
of $\sigma_2$. Moreover, a branch of $\tau$ which
is incident on a switch in $\sigma_i$ and 
which is not contained in $\sigma_i$ is a small
branch contained in a once punctured monogon
component $C$ of $S_i-\sigma_i$, and it is the
unique branch of $\tau_i$ in this component. 

For $i=0,1,2$ 
choose a tight splitting
sequence $\{\sigma_i(j)\}_{0\leq j\leq s_i}$ 
issuing from $\sigma_i(j)=\sigma_i$ with the property that
there is a pseudo-Anosov element
$\phi_i$ of ${\cal M}(S_i)$ such that $\sigma_i(s_i)$
is shift equivalent to $\phi_i(\sigma_i)$.
The existence of such a splitting
sequence follows from the arguments in the
proof of Lemma \ref{lemma55}. We view
$\phi_i$ as a reducible element of 
${\cal M}(S)$ which can be represented
by a diffeomorphism of $S$ fixing 
$\gamma_i$ (or $\gamma_1$ and $\gamma_2$ for $i=0$) pointwise. 
Define $\tau_i=\phi_i(\tau)$; then
$\tau_i$ is \emph{carried} by $\tau$ and 
contains $\phi_i\sigma_i$ as a subtrack.
Moreover, there is a carrying map 
$\tau_i\to \tau$ which maps every branch
$b$ of $\phi_i\sigma_i<\tau_i$ \emph{onto}
$\sigma_i$ and which maps $\tau_i-\phi_i\sigma_i$
bijectively onto $\tau-\sigma_i$.
In the sequel we always assume that our carrying
maps have these properties.

For $k>0$ define $\zeta(k)=\phi_1\circ 
\phi_0^{2k}\circ\phi_2\circ\phi_0^{2k}$ (where $\circ$ means
composition, i.e. $a\circ b$ represents the mapping
class obtained by applying $b$ first followed by an
application of $a$).
We claim that for every $k>0$ the mapping class $\zeta(k)$ is 
pseudo-Anosov.
For this note first that $\zeta(k)\tau$ is carried
by $\tau$ for all $k$. Thus it follows as
in Lemma \ref{lemma53} that $\zeta(k)$ is pseudo-Anosov
if and only if $\cap_i\zeta(k)^i{\cal V}(\tau)$
consists of a single ray which fills up $S$.
By the considerations in the proof of Lemma 5.4, this
is the case if a carrying map 
$\zeta(k)^2\tau\to \tau$ maps every branch
of $\zeta(k)^2\tau$ \emph{onto} $\tau$ (compare
also \cite{K85}).

Now let $b$ be any branch of $\zeta(k)^2\tau$;
then $b\in \zeta^2(k)\sigma_i$ for $i=1$ or $i=2$. Assume 
first that $b\in \zeta(k)^2\sigma_2$.
Since the splitting sequence $\{\sigma_2(j)\}_{0\leq i\leq s_2}$
is tight by assumption,
the image of $b$ under a
carrying map $\zeta(k)^2\tau\to \zeta(k)\phi_1\phi_0^{2k}\tau$ 
equals the subtrack $\zeta(k)\phi_1\phi_0^{2k}\sigma_2$
of $\zeta(k)\phi_1\phi_0^{2k}\tau$. In particular,
it contains the subgraph $\zeta(k)\phi_1\phi_0^{2k}(\tau-\sigma_1)$.
On the other hand, for the same reason
every branch of 
$\zeta(k)\phi_1\phi_0^{2k}\sigma_0<\zeta(k)\phi_1\phi_0^{2k}\sigma_2$ 
is mapped by a
carrying map $\zeta(k)\phi_1\phi_0^{2k}\tau\to
\zeta(k)\tau$ \emph{onto} $\zeta(k)\sigma_1$, and its
maps $\zeta(k)\phi_1\phi_0^{2k}(\tau-\sigma_1)$ 
bijectively onto $\zeta(k)(\tau-\sigma_1)$. It follows
that a carrying
map $\zeta(k)^2\tau\to \zeta(k)\tau$ maps $b$ 
\emph{onto} $\zeta(k)\tau$. 
Now if $b\in \zeta(k)^2\sigma_1$ then
the same argument yields that the image of $b$
under a carrying map $\zeta(k)^2\tau
\to \zeta(k)\tau$ contains the subtrack
$\zeta(k)\sigma_1$ of $\zeta(k)\tau$. 
In particular, it contains
a branch $c\in \zeta(k)\sigma_0< \zeta(k)\sigma_2$
which is mapped by a carrying map $\zeta(k)\tau
\to \tau$ \emph{onto} $\tau$. As before, we conclude
that a carrying map $\zeta(k)^2\tau\to \tau$
maps every branch $b$ of $\zeta(k)^2\tau$ onto $\tau$
and consequently each of the maps $\zeta(k)$
is pseudo-Anosov.

To establish our proposition it is enough to show
that for every
$\epsilon >0$ there is a number $k(\epsilon)>0$
with the following property.
For every $k\geq k(\epsilon)$,
the periodic orbit for the
Teichm\"uller geodesic flow on ${\cal Q}(S)$ 
which corresponds to the conjugacy class
of $\zeta(k)$ is entirely contained in the
subset of  
${\cal Q}(S)$ of all quadratic differentials
$q$ which admit an essential simple closed curve
of $q$-length at most $\epsilon$
(compare the proof of Lemma \ref{lemma55}).

By Lemma \ref{lemma52}, 
since for $i=1,2$ the splitting sequence 
$\{\sigma_i(j)\}_{0\leq j\leq s_i}$ is tight
there is a number $c>0$ with the following 
property. Let $\mu$ be a measured geodesic
lamination on $S_i$ which is carried by $\sigma_i(s_i)$
and which defines the transverse measure
$\mu_0\in {\cal V}(\sigma_i(0))$;
then $\mu_0(b_1)/\mu_0(b_2)\leq c$ for 
any two branches $b_1,b_2$ of $\sigma_i$.
Now the preimage of $\sigma_i$ under any
carrying map $\phi_i\tau\to \tau$ equals the
subtrack $\phi_i\sigma_i$ of $\phi_i\tau$ 
and the restriction
of a suitably chosen carrying map to $\phi_i\tau-\phi_i\sigma_i$
is injective. This implies that for
\emph{every} measured geodesic lamination
$\mu$ which is carried by $\phi_i\tau$ and
defines a transverse measure on $\phi_i\tau$ which is
not supported in $\phi_i\tau-\phi_i\sigma_i$, 
the measure $\mu_0$ on $\tau$ induced
from $\mu$ by a carrying map 
$\phi_i\tau\to \tau$ satisfies 
$\mu_0(b_1)/\mu_0(b_2)\leq c$ for all $b_1,b_2\in 
\sigma_i$.

For $\epsilon >0$ and $i=1,2$ let $C_i(\epsilon)$
be the closed subset of  ${\cal V}(\tau)$
containing all transverse measures $\nu$
with the following properties.
\begin{enumerate}
\item
The total weight of $\nu$ is one and the
sum of the $\nu$-weights over all 
branches of $\tau$ which are \emph{not}
contained in $\sigma_i$ is at most $\epsilon$.
\item For any two branches $b_1,b_2$ of 
$\sigma_i<\tau$ we
have $\nu(b_1)/\nu(b_2)\leq c$.
\end{enumerate}
For a transverse measure $\mu$ on $\tau$
denote as before by $\omega(\tau,\mu)$ the
total weight of $\mu$.
We claim that for every $\epsilon >0$
the transformation $\rho_i^k:{\cal V}(\tau)-\{0\}\to
{\cal V}(\tau)$ which is the composition
of the map $\phi_i\phi_0^{2k}:{\cal V}(\tau)\to 
{\cal V}(\phi_i\phi_0^{2k}\tau)$ with the map
${\cal V}(\phi_i\phi_0^{2k}\tau)\to 
{\cal V}(\tau)$ induced by a
carrying map $\phi_i\phi_0^k\tau\to \tau$
and with the normalization map $\mu\to \mu/\omega(\tau, \mu)$
maps $C_{i+1}(\epsilon)$ into to $C_i(\epsilon)$
provided that $k$ is sufficiently large
(indices are taken modulo 2).

We show our claim for the map $\rho_1^k$, the
claim for $\rho_2^k$ follows in exactly
the same way. Note first that by the choice
of $c$, the image under $\rho_1^k$ of \emph{every}
transverse measure on $\tau$ which
does not vanish on $\sigma_1$
satisfies property 2) above for $i=1$. Thus since
a measure from the set $C_2(\epsilon)$
does not vanish on $\sigma_0<\sigma_1$, property
2) holds for all measures in
$\rho_1^kC_{2}(\epsilon)$ and all $k>0$.
To establish property 1) for sufficiently large $k$, let
$\nu\in C_2(\epsilon)$, viewed as a measured geodesic
lamination, let $k>0$ and 
let $\mu=\phi_1\phi_0^{2k}(\nu)$. For
$0\leq i\leq 2k$ let $\mu(i)$ be the
transverse measure on 
$\phi_1\phi_0^i(\tau)$ of total weight one
which defines a multiple of the measured geodesic
lamination $\mu$.
Then there is some $a(i)\leq 1$ such that 
$\mu(i)=a(i)\mu$ as measured geodesic laminations.
By the definition of 
the set $C_2(\epsilon)$, there is a number
$r>0$ only depending on $\epsilon$ and $c$ 
such that the total $\mu$-weight of the subtrack
$\phi_1\phi_0^{2k}(\sigma_0)$ of $\phi_1\phi_0^{2k}(\tau)$ 
is not smaller than $r$.
Now the splitting sequence $\{\sigma_0(j)\}_{0\leq j\leq s_0}$
is tight and therefore a carrying map $\sigma_0(s_0)\to
\sigma_0(0)$ strictly increases the total weight of a 
transverse measure by at least a factor $L_0>1$.
Hence for every $s\leq 2k$
the total $\mu$-weight of 
$\phi_1\phi_0^{2k-s}(\sigma_0)<\phi_1\phi_0^{2k-s}(\tau)$ 
is not smaller than $1-r+L_0^sr$. Moreover, the sum
of the $\mu$-weights of the branches of 
$\phi_1\phi_0^{2k-s}(\tau)$ which are 
not contained in $\phi_1\phi_0^{2k-s}(\sigma_0)$ 
is independent of $s\leq 2k$. Thus the sum of the
$\mu(2k-s)$-weights of the branches
of $\phi_1\phi_0^{2k-s}(\tau)$ which are \emph{not} contained in 
$\phi_1\phi_0^{2k-s}(\sigma_0)$ is
smaller than $(1-r)/(1-r+L_0^sr)$. As a consequence,
if $(1-r)/(1-r+L_0^{k}r)<\epsilon$
then the $\mu(k)$-weight of
$\phi_1\phi_0^k(\tau-\sigma_0)$ does not exceed $\epsilon$
and the same holds true for the 
$\rho_1^k(\nu)$-weight of $\tau-\sigma_1$ which shows 
our above claim.

Together we deduce the existence of 
a number $k(\epsilon)>0$ such that
for $k\geq k(\epsilon)$ the set $C_1(\epsilon)$ is invariant under
the map $\rho$ which assigns to a measured 
geodesic lamination $0\not=\mu\in {\cal V}(\tau)$
the normalized
image of $\zeta(k)\mu\in {\cal V}(\zeta(k)\tau)$
under a carrying map ${\cal V}(\zeta(k)\tau)\to
{\cal V}(\tau)$. Now 
if $k>k(\epsilon)$ and if $(\lambda_k,\nu_k)\in {\cal B}(\tau)$
is such that
$q(\lambda_k,\nu_k)=q_k$ is a cotangent vector 
of the axis of $\zeta(k)$ then 
$\lambda_k$ spans the ray $\cap_i\zeta(k)^i{\cal V}(\tau)$ 
and therefore necessarily 
$\lambda_k\in C_1(\epsilon)$. 

Let $k\geq k(\epsilon)$ and 
let $a_k(1)>1$ be such that the total weight of the
transverse measure on $\phi_1\phi_0^k(\tau)$ defined
by the measured geodesic lamination
$a_k(1)\lambda_k$ equals one. Recall from Lemma 2.5
of \cite{H04} that
the intersection of $a_k(1)\lambda_k$ with
the embedded trainpath $\gamma_1$ on $\phi_1\phi_0^k(\tau)$
is bounded from above by the sum of
the $a_k(1)\lambda_k$-weights of the branches of
$\phi_1\phi_0^k(\tau)-\phi_1\phi_0^k(\sigma_1)$.
By our above consideration, we may assume that for
$k\geq k(\epsilon)$  this weight is bounded from above
by $\epsilon$. As a consequence, 
for every $t\leq \log a_k(1)$ 
the intersection of $\gamma_1$ with
the horizontal measured geodesic
lamination of the quadratic differential
$\Phi^tq_k$ is at most
$\epsilon$. 

Similarly, let $a_k(2)<a_k(3)<a_k(4)<a_k(5)$ be such that
the total weight of $a_k(2)\lambda_k$ on
$\phi_1\phi_0^{2k}(\tau)$ equals one, that the total weight of
$a_k(3)\lambda_k$ on  
$\phi_1\phi_0^{2k}\phi_2(\tau)$ equals one,
that the total weight of $a_k(4)\lambda_k$
on $\phi_1\phi_0^{2k}\phi_2\phi_0^k(\tau)$ equals
one and that the total weight
of $a_k(5)\lambda_k$ on $\zeta(k)\tau$ equals
one. Note that $\log a_k(5)$ is the length
of the periodic orbit of the Teichm\"uller flow
defined by the conjugacy class of $\zeta(k)$.
Using once more our above consideration
we conclude that for $k\geq k(\epsilon)$
and every $t\leq \log a_k(4)$ the intersection
of the horizontal measured geodesic lamination
of the quadratic differential 
$\Phi^tq_k$ with the curve $\zeta(k)\gamma_2$
is bounded from above by $\epsilon$.

Our proposition now follows if for sufficiently
large $k$ we can control
the intersections of the curves $\gamma_1,\zeta(k)\gamma_2$ 
with the vertical measured geodesic laminations
for the quadratic differentials $\Phi^tq_k$.
For this let again 
$\nu_k\in {\cal V}^*(\tau)$
be such that $(\lambda_k,\nu_k)\in {\cal B}(\tau)$
and $q_k=q(\lambda_k,\nu_k)$. 
Since a carrying map $\phi_1\phi_0^{2k}\phi_2(\tau)
\to \phi_1\phi_0^{2k}(\tau)$ maps every
branch $b$ of $\phi_1\phi_0^{2k}\phi_2(\sigma_2)$
\emph{onto} $\phi_0\phi_0^{2k}(\sigma_2)$ and
since $a_k(2)\lambda_k\in \phi_1\phi_0^{2k}C_2(\epsilon)$
by our above consideration, 
the $a_k(2)\lambda_k$-weight of 
every branch of $\phi_1\phi_0^{2k}(\tau)$ which
is contained in $\phi_1\phi_0^{2k}(\sigma_2)$
is bounded from below by a universal constant
$\delta>0$ not depending on $k$. 
Then the $a_k(2)^{-1}\nu_k$-weight of
every such branch is bounded from 
above by $1/\delta$ and hence there is a number
$\chi>0$ not depending on $k$ which
bounds from above the $a_k(3)^{-1}\nu_k$-weight of every
branch of $\phi_1\phi_0^{2k}\phi_2(\sigma_2)<
\phi_1\phi_0^{2k}\phi_2(\tau)$. Now the
curve $\phi_1\phi_0^{2k}\phi_2(\gamma_1)$ is 
an embedded trainpath in the train track 
$\phi_1\phi_0^{2k}\phi_2(\sigma_2)<
\phi_1\phi_0^{2k}\phi_2(\tau)$ 
and hence our upper bound for the
values of $a_k(3)^{-1}\nu_k$ on the branches of 
$\phi_1\phi_0^{2k}\phi_2(\sigma_2)$ implies that
the intersection between
$\phi_1\phi_0^{2k}\phi_2(\gamma_1)$ and 
$a_k(3)^{-1}\nu_k$ is uniformly bounded.
As we increase $k$, the ratios $a_k(4)/a_k(3)$ tend
to infinity (compare the above consideration) and hence
after possibly increasing
$k(\epsilon)$ we may assume that for 
every $k\geq k(\epsilon)$
we have $i(a_k(4)^{-1}\nu_k,\phi_1\phi_0^{2k}\phi_2\phi_0^{k}(\gamma_1))
<\epsilon$. By invariance of the intersection
form under the action of ${\cal M}(S)$ we conclude that
for every $k\geq k(\epsilon)$ and every
$t\geq \log(a_k(4)/a_k(5))$ the intersection between the
vertical measured geodesic lamination
of $\Phi^tq_k$ and $\gamma_1$ is bounded from above
by $\epsilon$. This then shows that for every
$t\in [\log(a_k(4)/a_k(5)),\log a_k(1)]$
the sum of the intersection numbers between $\gamma_1$
and the vertical and the horizontal measured 
geodesic lamination of $\Phi^t q_k$ does not exceed $2\epsilon$.
In other words, for every such $t$ the $\Phi^tq_k$-length
of $\gamma_1$ is bounded from above by $4\epsilon$.

The same argument shows that after possibly
increasing $k(\epsilon)$ once more we may assume that
for $k\geq k(\epsilon)$
and every $t\in [\log a_k(1),\log a_k(4)]$ the 
$\Phi^tq_k$-length of $\zeta(k)(\gamma_2)$ is bounded from
above by $4\epsilon$. By periodicity, we conclude that
for $k\geq k(\epsilon)$ 
the periodic orbit of $\Phi^t$ which corresponds to
the conjugacy class of $\zeta(k)$ 
is entirely contained in the set of quadratic
differentials $q$ which admits an essential simple
closed curve of $q$-length at most $4\epsilon$.
This completes the proof of our proposition.
\end{proof}

\bigskip

{\bf Remark:} 1) A pseudo-Anosov element $\phi$
acts as an isometry on the \emph{curve graph} 
$({\cal C}(S),d)$ of $S$. By a result of Bowditch
\cite{Bw03}, for every $c\in {\cal C}(S)$ the limit 
$\lim_{k\to \infty}d(\phi^kc,c)/k$ exists 
and is independent of $c$. This limit is called
the \emph{stable length} for this action.
The stable length of each of the (infinitely many)
pseudo-Anosov elements $\zeta(k)$ constructed
in the proof of Proposition 5.6 is at most 2.
Moreover, let $\gamma_k$ be the 
Teichm\"uller geodesic
in ${\cal T}(S)$ which is invariant
under $\zeta(k)$. Then there is a number
$\epsilon >0$ such that 
for sufficiently large $k$, the set of
essential 
simple closed curves $c$ on our surface $S$ for which
the minimum of the hyperbolic lengths 
of $c$ along $\gamma_k$ is smaller than $\epsilon$
is precisely an orbit under the action of
$\zeta(k)$ of a pair of disjoint
essential simple closed curves on $S$.

2) We believe that Proposition \ref{prop56} is valid for
\emph{every} nonexceptional surface of finite type, 
with pseudo-Anosov mapping classes which
can be constructed as the once in the proof our proposition. 
However 
we did not attempt to carry out the details.

\bigskip

\noindent
MATHEMATISCHES INSTITUT DER UNIVERSIT\"AT BONN\\
BERINGSTRASSE 1, D-53115 BONN, GERMANY

\noindent
e-mail: ursula@math.uni-bonn.de


\begin{thebibliography}{CEG87}










\bibitem[Bw03]{Bw03} B.~Bowditch, {\em Tight
geodesics in the curve complex}, Invetn. Math. 171 (2008),
281--300.



\bibitem[CEG87]{CEG87} R.~Canary, D.~Epstein, P.~Green,
{\em Notes on notes of Thurston}, in ``Analytical and geometric
aspects of hyperbolic space'', edited by D.~Epstein, London Math.
Soc. Lecture Notes 111, Cambridge University Press, Cambridge 1987.

\bibitem[CB88]{CB88} A.~Casson with S.~Bleiler, {\sl Automorphisms
of surfaces after Nielsen and Thurston}, Cambridge University
Press, Cambridge 1988.



\bibitem[FLP91]{FLP91} A.~Fathi, F.~Laudenbach, V.~Po\'enaru, {\sl Travaux de
Thurston sur les surfaces,} Ast\'erisque 1991.

\bibitem[H06]{H04} U.~Hamenst\"adt, {\em Train
tracks and the Gromov boundary of the complex
of curves}, in ``Spaces of Kleinian groups''
(Y.~Minsky, M.~Sakuma, C.~Series, eds.),
London Math. Soc. Lec. Notes 329 (2006), 187--207,
Cambridge University Press, Cambridge 2006.

\bibitem[H09]{H05} U.~Hamenst\"adt, {\em Geometry of the
mapping class groups I: Boundary amenability}, 
Invent. Math. 175 (2009), 545--609.







\bibitem[I02]{I02} N.~V.~Ivanov, {\it Mapping class groups},
Chapter 12 in Handbook of Geometric Topology (Editors
R.J.~Daverman and R.B.~Sher), Elsevier Science (2002), 523-633.

\bibitem[K85]{K85} S.~Kerckhoff, {\em Simplicial
systems for interval exchange maps and measured foliations},
Erg. Th. \& Dyn. Sys. 5 (1985), 257--271.

\bibitem[K92]{K92} S.~Kerckhoff, {\em Lines of
minima in Teichm\"uller space}, Duke Math. J. 65 (1992), 187--213.




\bibitem[L83]{L83} G.~Levitt, {\em Foliations and
laminations on hyperbolic surfaces}, Topology 22 (1983),
119--135.





\bibitem[MM99]{MM99} H.~Masur, Y.~Minsky, {\em Geometry of the
complex of curves I: Hyperbolicity}, Invent. Math. 138 (1999),
103-149.





\bibitem[Mo03]{M03} L.~Mosher, {\em Train track expansions of measured
foliations}, unpublished manuscript.

\bibitem[O96]{O96} J.~P.~Otal, {\sl Le Th\'{e}or\`{e}me d'hyperbolisation
pour les vari\'et\'es fibr\'ees de dimension 3}, Ast\'erisque 235,
Soc. Math. Fr. 1996.

\bibitem[PH92]{PH92} R.~Penner with J.~Harer, {\sl Combinatorics
of train tracks}, Ann. Math. Studies 125, Princeton University
Press, Princeton 1992.

\bibitem[R05]{R05} K.~Rafi, {\em A characterization of short
curves of a Teichm\"uller geodesic}, Geom. Top. 9 (2005),
179--202.







\end{thebibliography}
\end{document}